# Applications of the Digital-Discrete Method in Smooth-Continuous Data Reconstruction

Li Chen
Department of Computer Science and Information Technology
University of the District of Columbia
Email: *lchen@udc.edu*

**Abstract:** This paper presents some applications of using recently developed algorithms for smooth-continuous data reconstruction based on the digital-discrete method. The classical discrete method for data reconstruction is based on domain decomposition according to guiding (or sample) points. Then uses Splines (for polynomial) or finite elements method (for PDE) to fit the data. Our method is based on the gradually varied function that does not assume the property of the linearly separable among guiding points, i.e. no domain decomposition methods are needed. We also demonstrate the flexibility of the new method and the potential to solve variety of problems. The examples include some real data from water well logs and harmonic functions on closed 2D manifolds. This paper presented the results from six different algorithms. This method can be easily extended to higher multi-dimensions.

## 1. Introduction

Data reconstruction is to fit a function based on the observations of some sample (guiding) points. This paper presents some applications of using recently developed algorithms for smooth-continuous data reconstruction based on the digital-discrete method. The classical discrete method for data reconstruction is based on domain decomposition according to guiding (or sample) points. Then uses Splines (for polynomial) or finite elements method (for PDE) to fit the data.

Some successful methods have been discovered or proposed to solve the problem including the moving least square method and others presented in [1][10][13][14].

Our method is based on the gradually varied function that does not assume the property of the linearly separable among guiding points, i.e. no domain decomposition methods are needed. We also demonstrate the flexibility of the new method and the potential to solve variety of problems. The examples include some real data from water well logs and harmonic functions on closed 2D manifolds. This paper presented the results from six different algorithms. This method can be easily extended to higher multi-dimensions.

This paper will present applications of using the method to many different cases including the rectangle domain and closed surfaces. In [2], we directly deal with smooth functions on the (same) piecewise linear manifold or a non-Jordon graph. We also have applied this method to groundwater flow equations [3]. We have discovered that the major difference between our methods to existing methods is that the former is a true nonlinear approach [17]. We also have analyzed the relationships and differences among different method in [17].

## 2. Basic Concepts

Gradual variation is a discrete method that can be built on any graph. The gradually varied surface is a special discrete surface. We now introduce this concept.

The Concept of Gradual Variation: Let function f: D→{A1, A2,…,An}, if a and b are adjacent in D implies f(a)=f(b), or f(b) =A(i-1) orA(i+1) when f(a)=Ai , point (a,f(a)) and (b,f(b)) are said to be gradually varied. A 2D function (surface) is said to be gradually varied if every adjacent pair are gradually varied.

Discrete Surface Fitting: Given J⊆D, and f: J→{A1,A2,…An} decide if there is a F: D→{A1,A2,…,An} such that F is gradually varied where f(x)=F(x), x in J.

**Theorem** (Chen, 1989) The necessary and sufficient conditions for the existence of a gradually varied extension F is: for all x,y in J, d(x,y)≥ |i-j|, f(x)=Ai and f(y)=Aj, where d is the distance between x and y in D.

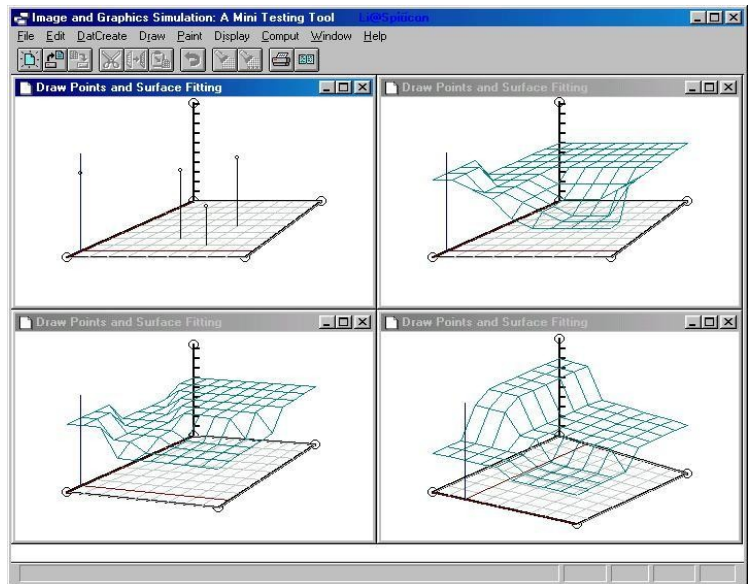

Fig. 2 Examples of gradually varied functions

A gradually varied surface fitting software component was included in a lab-use oriented software system in 1997. The software can demo the arbitrary guiding points gradually varied surface fitting in a 10x10 grid domain.

In 2004, Chen and Adjei proposed a method to do continuous and differential surface fitting based on lambda-connectedness that is a continuous space treatment of discrete functions [8].

However, this method has not been implemented. In 2008, Chen implemented an algorithm for groundwater data reconstruction; the implementation was not fully successful.

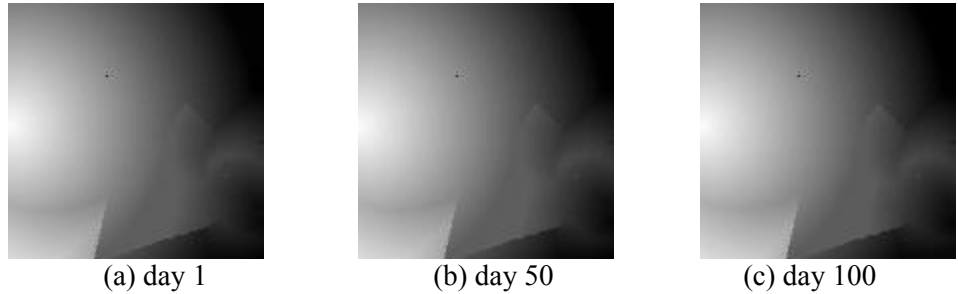

    (a) day 1        (b) day 50       (c) day 100

Fig. 3. Northern VA Groundwater distribution calculated by gradually varied surfaces date from 04/01/07. The intensity indicated the depth of the groundwater. 10 sample points are used.

Thus, the above theorem can be used for a single surface fitting if the condition in the theorem is satisfied. The problem is that the sample data does not satisfy the condition of fitting. So the original algorithm cannot be used directly for individual surface fitting. Another problem is that the theorem is only for "continuous" surfaces. It does not imply a solution for differtiable or smooth functions.

## 3. Algorithms and Experiments

In [2], a systematic digital-discrete method for obtaining continuous functions with smoothness to a certain order ($C^n$) from sample data is designed. This method is based on gradually varied functions and the classical finite difference method. This new method has been applied to real groundwater data and the results have validated the method. This method is independent from existing popular methods such as the cubic spline method and the finite element method. The new digital-discrete method has considerable advantages for a large amount of real data applications. This digital method also differs from other classical discrete method that usually uses triangulations. This method can potentially be used to obtain smooth functions such as polynomials through its derivatives $f^{(k)}$ and the solution for partial differential equations such as harmonic and other important equations.

The new algorithm tries to search for a best solution of the fitting. We have added the component of the classical finite difference method. The major steps of the new algorithm include: (This is for 2D functions. For 3D function, we only need to add a dimension).

Step 1: Load guiding points. In this step we load the data points with observation values.
Step 2: Determine resolution. Locate the points in grid space.
Step 3: Function extension according to the theorem presented in Section 2. So we get the
    gradually varied or near gradually varied (continuous) functions.
Step 4: Use finite difference method to calculate partial derivatives. Then get the smoothed
    function.
Step 5: Some multilevel and multi resolution method may be used.

Three sets of real data are tested. The second set used the same 10 points sample, the result is improved.

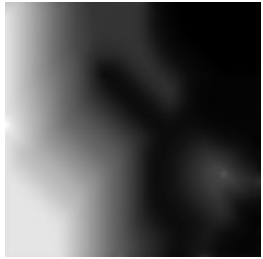

Fig. 4. Northern VA Groundwater distribution calculated by gradually varied surfaces date from 04/01/07. 10 sample points are used at Day 95.

The following example shows the third data set containing 29 sample points.

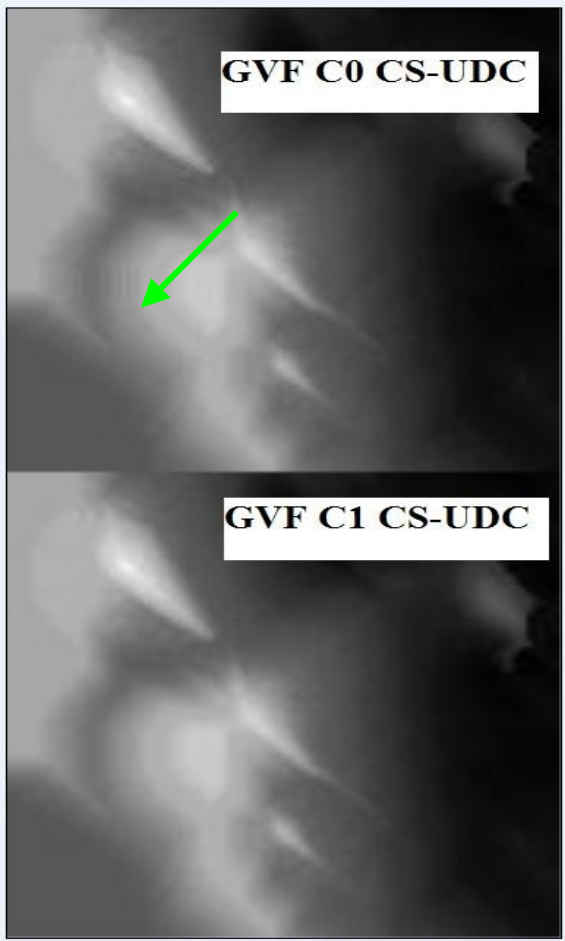

Fig.5. The picture is the result of fitting based on 29 sample points. The first is a "continuous" surface and the second one has "first derivative." The arrow indicated the interesting area disappeared in second image.

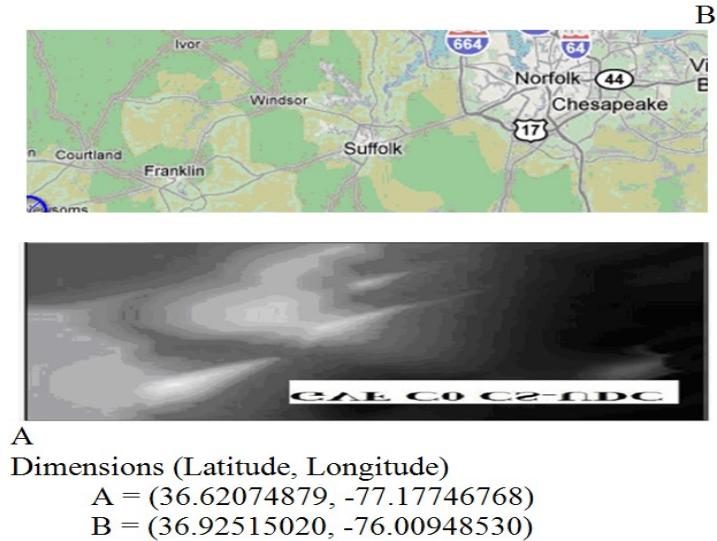

Fig. 6. The map and ground water data

Fig.6. It shows a good match found between the ground water data and the region geographical map. The brighter pixels mean the deeper distance from the surface. In mountain area, the groundwater level is lower in general. Some mismatches may be caused by not having enough sample data points (wells).

## 4. Continuous and Smooth Functions on Manifolds

A gradually varied surface reconstruction does not rely on the shape of domain. And it does not be restricted by simplicial decomposition. As long as the domain is described as a graph, our algorithms will apply. However, the actual implementation will be much difficult. In the above sections, we have discussed two types of algorithms for a rectangle domain. One is for the complete gradually varied function (GVF) fitting, and another is to reconstruct the best fit based on the gradual variation and the finite difference method.

The following is the implementation of the method for digital-discrete surface fitting on manifolds (triangulated representation for the domain). Data came from a modified example in Princeton 3D Benchmark data sets.

We will have four algorithms related to continuous (and smooth) functions on manifolds. This is because that we have 4 cases: (1) ManifoldIntGVF: The GVF extension on point space, corresponding to Delaunay triangulations; the values are integers. (2) ManifoldRealGVF: The GVF fitting on point space, the fitted data are real numbers. (3) ManifoldCellIntGVF: The GVF extension on face (2D-cell) space, corresponding to Voronoi decomposition; the values are integers. (4) ManifoldCellRealGVF: The GVF fitting on face (2D-cell) space, the fitted data are real numbers.

Note that ManifoldRealGVF and ManifoldCellRealGVF are the algorithms that based on gradually varied functions. The results may or may not be gradually varied. A special data structure was chosen to hold all points, edges, and 2D-cells. It is not very easy to demonstrate the correctness of the reconstruction on manifolds since to determine the location on 3D display is hard.

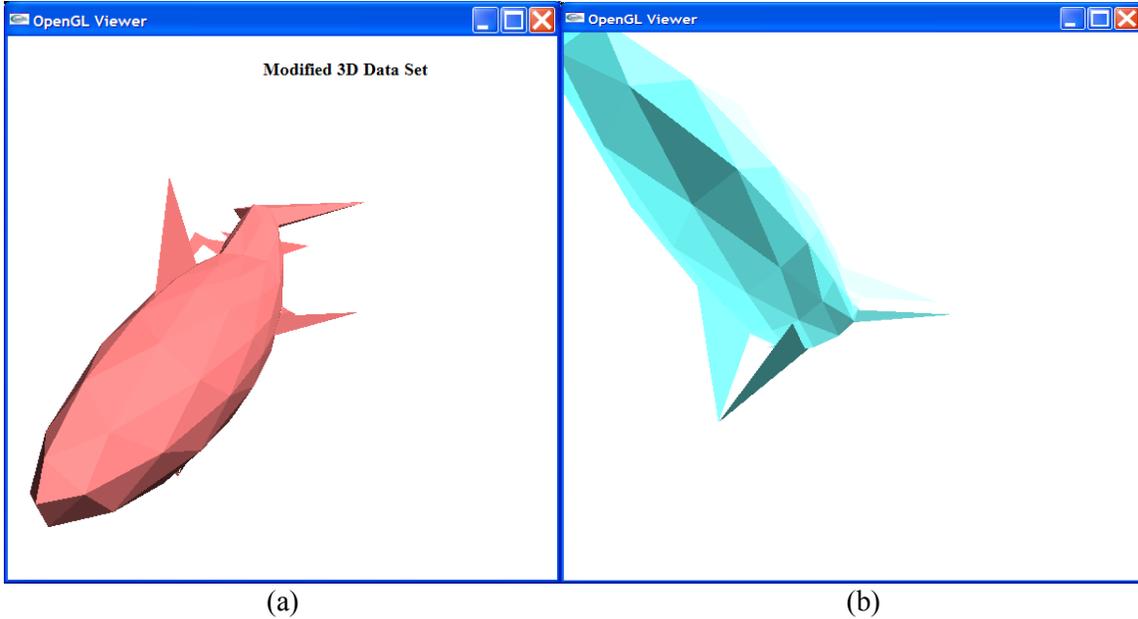

(a)                                                                                    (b)

Fig. 7 The ManifoldIntGVF Algorithm: (a) The left is the original 3D image, (b) the reconstruction based on four sample points.

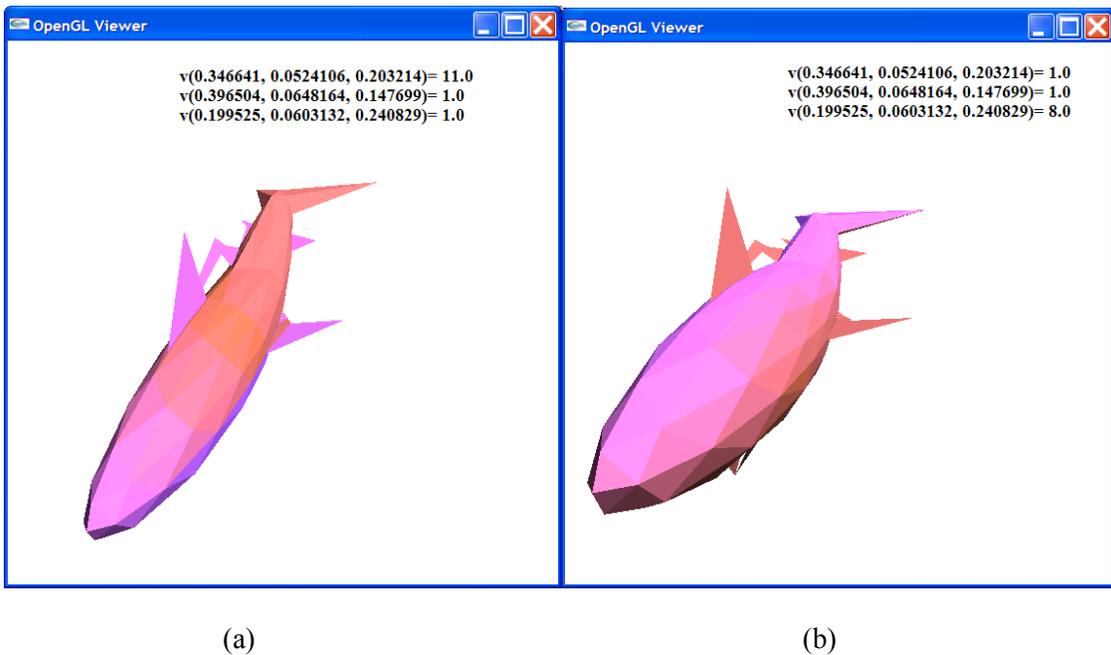

(a)                                                                                    (b)

Fig. 8. The ManifoldIntGVF Algorithm for two 3 sample points sets: (a) The gradual variation is from back to front, (b) The gradual variation from top to right side.

In Fig. 7 (a), we show an example using ManifoldIntGVF about 2D closed manifold (can be view as the boundary of a 3D object). It is a fish from Princeton 3D data sets. We have modified the data by adding some triangles to make the graph to be connected. Then we put four values on four vertices. Fig. 7 (b) shows the result. In order to avoid that we may be fooled because the 3D display might cause the color change, we have used other two sets of guiding points. Each of these uses only three guiding points. The second set of the guiding points is just swapped a pair of values with the first set. The values are posted in the pictures. See Fig. 8.

It is actually hard to select values for a relatively large set of sample points for a complete GVF fitting such as ManifoldIntGVF to satisfy the condition of gradual variation (actually the Lipschitz condition). For ManifoldRealGVF, we can choose any data values. Fig. 9 shows the results with 6 guiding points.

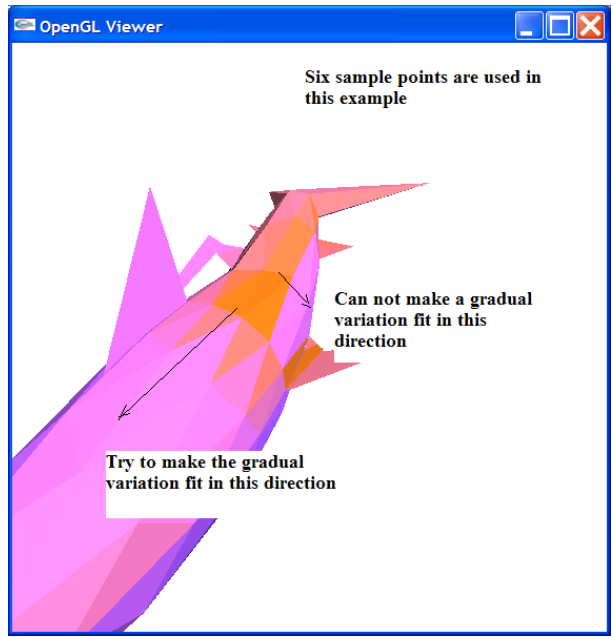

Fig. 9 The ManifoldRealGVF algorithm using 6 guiding points.

It is actually now we can get smooth functions based on the GVF fitting. We can use the Catmull-Clark method [10] or other methods. A small problem we got is the display. In order to make the display for the vertices, we use the average value of three vertices to be the display value on the triangle. We will make smaller triangle later to make more precise display later. This problem will disappear in the cell-based GVF. However, it needs to calculate the graph first for the cells and their adjacent cells.

The following examples are related to ManifoldCellIntGVF and ManifoldCellRealGVF. In order to make harmonic functions on manifolds, we restrict here the function to be gradually varied. Some discussions about harmonic functions and gradually varied function were presented in [9]. In the application of this paper, the maximum slop of two sample points are used to determine the value levels. This treatment will guarantee the sample points set satisfy the gradual variation condition.

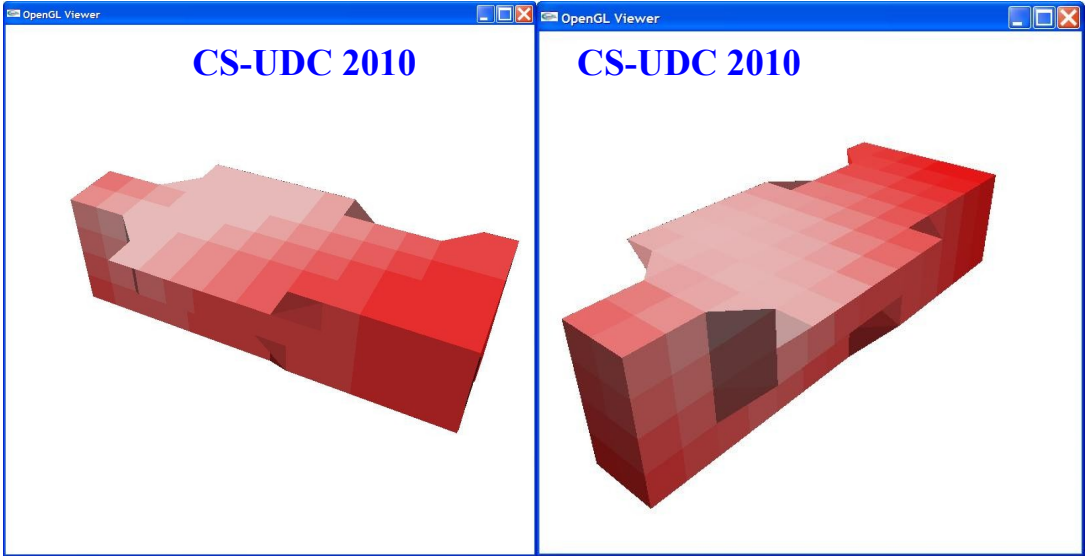

(a)                            (b)

Fig. 10. Using 7 Points to fit the data on 3D surface: (a) the GVF result. (b) The Harmonic fitting (iterations few times) based on GVF.

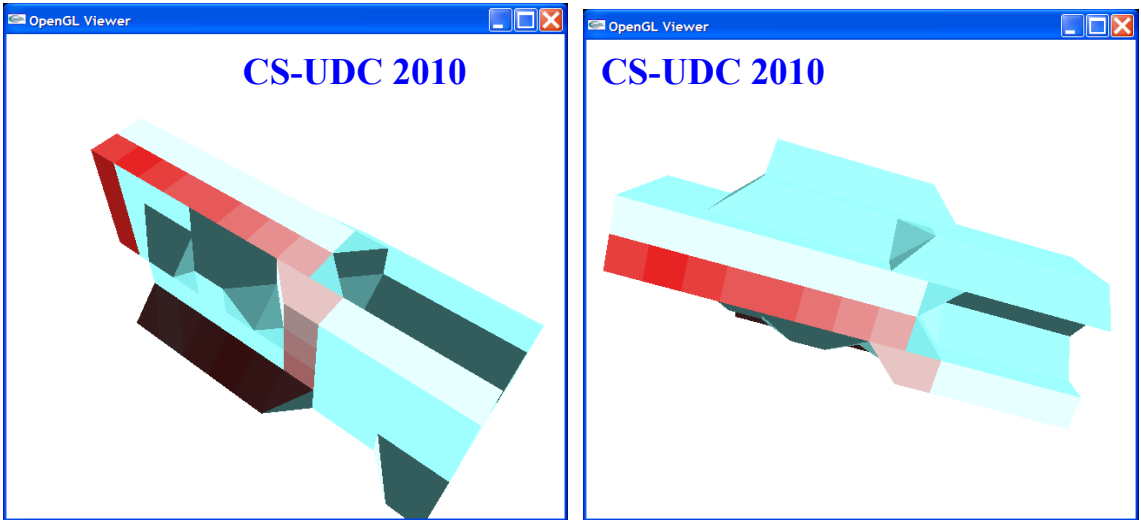

(a)                            (b)

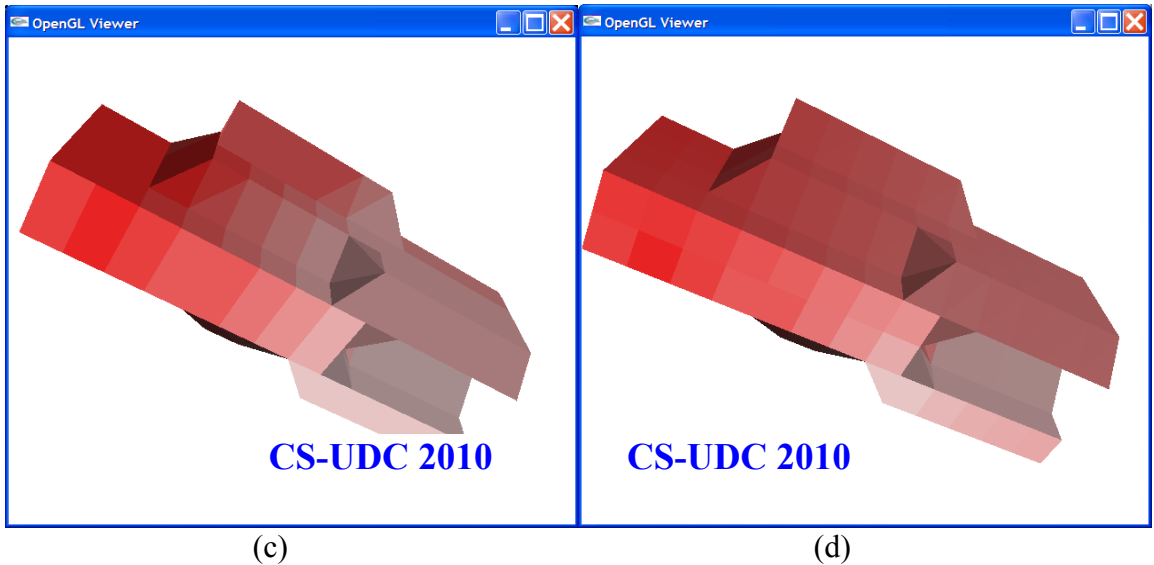

Fig. 11. The selected cells form a boundary curve that is gradually varied: (a) and (b) are two displays for the guiding points (cells). (c) The GVF result. (d) The Harmonic fitting based on GVF (100 iterations).

The smooth function on piece-wise linear manifolds is often studied in computer graphics. The main difference from our method and existing methods in graphics is that we do not know the value on the each vertex (or faces). Our method can solve this key problem of interpolation. After the values at all vertices are known for each vertex, we can construct the smooth functions for geodesic paths that cover all edges link to a vertex. Then use vector space calculus to get the unknown points. Partition of unity may be used here for getting a smooth solution on partition edges (similar to some finite element algorithms). We can also use the existing methods for obtaining the smooth functions [14][13].

In addition, we have used our algorithm to fit the data with 7 sample points on the manifold, as described earlier in Fig 10. We also use discrete harmonic functions to fit the data [9]. Since if boundary is known, then the harmonic solution will be unique, so we test an example for such a case in Fig. 11.

## 5. Summary

To get a smoothed function using gradual variation is a long time goal of our research. Some theoretical attempts have been made before, but struggled on the actual implementation. The author was invited to give a talk at the Workshop on the Whitney's Problem organized by Princeton University and College of William and Mary in Aug. 2009. He was somewhat inspired and encouraged by the presentations and the helpful discussions with the attendees of this workshop.

The purpose of this paper is to present some actual examples and related results using the new algorithms we designed in [2]. The author welcomes other real data sets to examine the new algorithms. The implementation code is written in C++. Li Chen's website is at www.udc.edu/prof/chen.

*Acknowledgement:* This research has been partially supported by the USGS seed grants through UDC Water Resources Research Institute (WRRI) and Center for Discrete Mathematics and Theoretical Computer Science (DIMACS) at Rutgers University. Professor Feng Luo suggested the direction of the relationship between harmonic functions and gradually varied functions. UDC undergraduate Travis Branham extracted the application data from the USGS database. Professor Thomas Funkhouser provided helps on the 3D data sets and OpenGL display programs. The author would also like to thank Professor C. Fefferman and Professor N. Zobin for their invitation for the Workshop on the Whitney's Problem in 2009.